\documentclass[10pt,a4paper]{amsart}
\usepackage[a4paper,marginratio=1:1]{geometry}
\usepackage[initials,non-sorted-cites]{amsrefs}
\usepackage{mybibspec}
\usepackage[symbol]{footmisc}
\usepackage{braket,amssymb}

\numberwithin{equation}{section}

\newcommand{\bdry}{\partial}
\newcommand{\abs}[1]{\lvert#1\rvert}
\newcommand{\CSp}{\mathit{CSp}}
\newcommand{\CU}{\mathit{CU}}
\DeclareMathOperator{\id}{id}
\DeclareMathOperator{\tr}{tr}
\DeclareMathOperator{\End}{End}
\DeclareMathOperator{\Ric}{Ric}

\theoremstyle{theorem}
\newtheorem{thm}{Theorem}[section]
\theoremstyle{definition}
\newtheorem{dfn}[thm]{Definition}

\title[Einstein metrics on strictly pseudoconvex domains]{Einstein metrics on strictly pseudoconvex domains from the viewpoint of bulk-boundary correspondence}
\author{Yoshihiko Matsumoto}
\address{Department of Mathematics, Graduate School of Science, Osaka University, Toyonaka, Osaka 560-0043, Japan}
\address{Department of Mathematics, Stanford University, Stanford, CA 94305-2125, USA}
\email{matsumoto@math.sci.osaka-u.ac.jp}

\begin{document}

\maketitle

\section{Introduction}
\label{sec:Introduction}

This is a survey discussing some aspects of the correspondence proposed
by Biquard \cites{Biquard-99,Biquard-00} (see also \cite{Biquard-00-English})
between Einstein metrics on the interior of
a manifold-with-boundary and strictly pseudoconvex CR structures on the boundary.
Here, by ``CR structures,'' we shall mean not only integrable almost CR structures
but also certain nonintegrable ones, which will be described later in this section.

We may regard this correspondence as a differential-geometric interpretation and generalization
of the classical complex-analytic correspondence
between strictly pseudoconvex domains and CR manifolds arising as their boundaries.
Let us start with describing this viewpoint.

Fefferman's mapping theorem \cite{Fefferman-74} states that
any biholomorphism $\Omega_1\to\Omega_2$ between smoothly bounded
strictly pseudoconvex domains in $\mathbb{C}^n$ ($n\ge 2$)
extends to a diffeomorphism between the closures of the domains.
This extension automatically restricts to a CR-diffeomorphism from $\bdry\Omega_1$ to $\bdry\Omega_2$.
Conversely, any CR-diffeomorphism $\bdry\Omega_1\to\bdry\Omega_2$ necessarily extends to
a biholomorphism $\Omega_1\to\Omega_2$ by the Bochner--Hartogs theorem \cite{Bochner-43}.
Such phenomena for domains in $\mathbb{C}^n$ generalize to those in Stein manifolds by the works of
Bedford--Bell--Catlin \cite{Bedford-Bell-Catlin-83} and Kohn--Rossi \cite{Kohn-Rossi-65},
and as a consequence, if $\mathcal{D}$ denotes the set of
all smoothly bounded strictly pseudoconvex domains in Stein manifolds,
then biholomorphism classes of domains in $\mathcal{D}$ and
CR-diffeomorphism classes of the boundaries of domains in $\mathcal{D}$ are in
a one-to-one correspondence:
\begin{equation}
	\label{eq:complex-CR-correspondence}
	{\mathcal{D}/\sim_{\text{bihol}}}\cong
	\set{\bdry\Omega|\Omega\in\mathcal{D}}/\sim_{\text{CR-diffeo}}.
\end{equation}

The classical approach toward the correspondence \eqref{eq:complex-CR-correspondence}
from differential geometry uses the Bergman metric.
However, here we would rather make use of the Einstein metric of Cheng--Yau \cite{Cheng-Yau-80}
in the following theorem,
because the Einstein equation has an advantage that it makes sense without complex structures
(recall that we are going to include some nonintegrable CR structures on the boundary into our consideration).

\begin{thm}[Cheng--Yau \cite{Cheng-Yau-80}]
	\label{thm:Cheng-Yau}
	Let $\Omega$ be a smoothly bounded strictly pseudoconvex domain in a Stein manifold of dimension $n\ge 2$.
	Then there exists a complete K\"ahler-Einstein metric with negative Einstein constant on $\Omega$,
	which is unique up to homothety.
\end{thm}

The metric is determined by the complex structure of $\Omega$.
On the other hand, the asymptotic behavior of the metric at the boundary
is mostly determined by the local CR geometry of $\bdry\Omega$,
as discussed by Fefferman \cite{Fefferman-76} and Graham \cite{Graham-87}.
In this sense, the Cheng--Yau metric is a link that realizes (to some extent) the correspondence
\eqref{eq:complex-CR-correspondence}. % in a quantitative manner.
We will introduce the notion of ``asymptotically complex hyperbolic Einstein metrics''
in the next section, for which the Cheng--Yau metrics serve as model examples.
%Our purpose is to establish an interesting class of complete Einstein metrics on non-compact manifolds,
%which is broader than that of the Cheng--Yau metrics,
%such that each metric is associated with some (possibly) non-integrable CR structure on the boundary at infinity.

Now we illustrate the class of almost CR structures that we consider.
Let $M$ be a (connected) differentiable manifold of dimension $2n-1$,
where $n\ge 2$, and $H$ a contact distribution over $M$.
We say, in this article, that an almost CR structure $J$ on $H$ (i.e., $J\in\Gamma(\End(H))$ satisfying
$J^2=-\id$) is \emph{compatible} when the \emph{Levi form}
\begin{equation*}
	h_\theta(X,Y):=d\theta(X,JY),\qquad X,\,Y\in H
\end{equation*}
is symmetric in $X$ and $Y$ for some (hence any) contact 1-form $\theta$ annihilating the distribution $H$.
(If $H^\perp\subset T^*M$ is oriented, then $H$ has a natural $\CSp(n-1)$-structure,
and fixing a compatible almost CR structure amounts to a reduction of the structure group of $H$
from $\CSp(n-1)$ to $\CU(p,n-1-p)$ for some $p$.
The term ``compatible'' refers to the compatibility of $J$ to the $\CSp(n-1)$-structure of $H$ in this setting.)
Integrable almost CR structures are always compatible as is well known,
but there are more compatible structures
(except for the three-dimensional case, in which any almost CR structure is automatically integrable).
It can be easily checked that $J$ is compatible if and only if
\begin{equation}
	\label{eq:partial-integrability}
	[\Gamma(H^{1,0}_J),\Gamma(H^{1,0}_J)]\subset\Gamma(H_\mathbb{C}),
\end{equation}
where $H_{\mathbb{C}}$ is the complexification of $H$ and
$H_{\mathbb{C}}=H^{1,0}_J\oplus H^{0,1}_J$ is the eigenbundle decomposition with respect to $J$
(note that \eqref{eq:partial-integrability} is not a trivial condition since
$H_\mathbb{C}\subsetneq T_\mathbb{C}\bdry\Omega$).
Because of \eqref{eq:partial-integrability},
compatible almost CR structures are also called \emph{partially integrable} in the literature
(e.g., \cite{Cap-Schichl-00,Cap-Slovak-09,Matsumoto-14,Matsumoto-16,Matsumoto-preprint16}).

The usual notion of strict pseudoconvexity naturally extends to compatible almost CR structures.
Namely, a compatible almost CR structure $J$
is said to be \emph{strictly pseudoconvex} if the Levi form $h_\theta$ has definite signature.
In the following we shall always assume the strict pseudoconvexity,
and a contact form $\theta$ is always taken so that $h_\theta$ is positive definite.

Each asymptotically complex hyperbolic (ACH for short) Einstein metric is ``associated'' to,
or ``fills inside'' of, a manifold equipped with a strictly pseudoconvex compatible almost CR structure,
as the Cheng--Yau metric does.
Relationships between Einstein metrics and geometric structures on the boundary have been more actively studied
in the setting of Poincar\'e-Einstein (or AH-Einstein) metrics and conformal structures,
partly because of physical interests.
Furthermore, the cases of Poincar\'e-Einstein and ACH-Einstein metrics
are generalized to a broader perspective involving ``asymptotically symmetric Einstein metrics'' and
``parabolic geometries,'' which is illustrated in \cites{Biquard-99,Biquard-00,Biquard-Mazzeo-06}.
The term ``bulk-boundary correspondence'' in the title of this article is intended to indicate this very general
correspondence, most part of which is yet to be unveiled.

\section{Asymptotically complex hyperbolic Einstein metrics}
\label{sec:ACH-metrics}

In order to motivate our definition of ACH metrics,
let us first observe the fact that the leading part of the asymptotic behavior of
the Cheng--Yau metric $g$ at the boundary can be described in terms of the CR structure of $\bdry\Omega$.

From the proof of its existence, it is known that $g$ is expressed (after a normalization) as
\begin{equation}
	\label{eq:Cheng-Yau-metric}
	g_{i\overline{j}}=\partial_i\partial_{\overline{j}}\log\frac{1}{\varphi}
	=\frac{(\partial_i\varphi)(\partial_{\overline{j}}\varphi)}{\varphi^2}
	-\frac{\partial_i\partial_{\overline{j}}\varphi}{\varphi},
\end{equation}
where $\varphi\in C^\infty(\Omega)\cap C^{n+1,\alpha}(\overline{\Omega})$ is some defining function
of $\Omega$, i.e.,
$\Omega=\set{\varphi>0}$ and $d\varphi$ is nowhere vanishing on $\bdry\Omega$,
where $\alpha\in(0,1)$ is arbitrary
(to be precise, \cite{Lee-Melrose-82} is responsible for this optimal boundary regularity).
Because of \eqref{eq:Cheng-Yau-metric}, one can take a diffeomorphism of the form
\begin{equation}
	\label{eq:collar-neighborhood-diffeomorphism}
	\Phi=(\pi,\rho)\colon\mathcal{U}\to\bdry\Omega\times[0,\varepsilon),
\end{equation}
where $\mathcal{U}$ is an open neighborhood of $\bdry\Omega$ in $\overline{\Omega}$,
such that $\pi\colon\mathcal{U}\to\bdry\Omega$ restricts to the identity map on $\bdry\Omega$ and
the Cheng--Yau metric $g$ satisfies
\begin{equation}
	\label{eq:asymptotic-complex-hyperbolicity}
	g\sim \Phi^*g_\theta,\qquad
	g_\theta=\frac{1}{2}\left(\frac{d\rho^2}{\rho^2}+\frac{\theta^2}{\rho^2}+\frac{h_\theta}{\rho}\right)
\end{equation}
as $\rho\to 0$,
where $\theta$ is some contact form on $\bdry\Omega$ annihilating the natural contact distribution
and $h_\theta$ is the associated Levi form.
The meaning of \eqref{eq:asymptotic-complex-hyperbolicity} can be understood as, for example, that
$\abs{g-\Phi^*g_\theta}_{\Phi^*g_\theta}$ uniformly tends to $0$ as $\rho\to 0$.
More is true actually: it follows from the asymptotic expansion established in \cite{Lee-Melrose-82}
that the $C^k$ norm of $\rho^{-1}(g-\Phi^*g_\theta)$,
defined geometrically by $\Phi^*g_\theta$, is finite for any $k\ge 0$.
We note that there exists, for any choice of $\theta$, a diffeomorphism $\Phi$ with respect to which
\eqref{eq:asymptotic-complex-hyperbolicity} holds---there is no preferred choice of $\theta$.
We also remark that in the literature
the model metric $g_\theta$ is sometimes expressed as
\begin{equation*}
	g_\theta=\frac{1}{2}\left(4\frac{dx^2}{x^2}+\frac{\theta^2}{x^4}+\frac{h_\theta}{x^2}\right)
\end{equation*}
by introducing a new coordinate $x=\sqrt{\rho}$.

Observing the asymptotic behavior \eqref{eq:asymptotic-complex-hyperbolicity} of the Cheng--Yau metric,
we define as follows.
Metrics of this type are firstly considered by Epstein--Melrose--Mendoza \cite{Epstein-Melrose-Mendoza-91},
in which the meromorphic continuation of the resolvent of the Laplacian (on functions) is studied.

\begin{dfn}
	Let $\overline{X}$ be a compact smooth manifold-with-boundary with $\dim_\mathbb{R}\overline{X}=2n$,
	where $n\ge 2$, and $X$ its interior.
	A Riemannian metric $g$ defined on $X$
	is called an \emph{asymptotically complex hyperbolic metric} (or \emph{ACH metric})
	when there exists a diffeomorphism $\Phi$ like \eqref{eq:collar-neighborhood-diffeomorphism}
	such that $g$ satisfies
	\eqref{eq:asymptotic-complex-hyperbolicity} with respect to some contact distribution $H$ over $\bdry X$,
	a strictly pseudoconvex compatible almost CR structure $J$ on $H$,
	and a contact form $\theta$ in the sense that
	\begin{equation*}
		g-\Phi^*g_\theta\in C^{2,\alpha}_\delta(X,S^2T^*X)
	\end{equation*}
	for some $\delta>0$ and arbitrary $\alpha\in(0,1)$.
	Here $C^{k,\alpha}_\delta(X,S^2T^*X)$ denotes the space of $C^k$ symmetric 2-tensors $\sigma$ on $X$
	such that $\rho^{-\delta/2}\sigma$ has finite $C^{k,\alpha}$ norm with respect to $\Phi^*g_\theta$
	(this space depends on $\Phi$ and $H$, but not on $J$).
	The almost CR structure $J$, or the triple $(\bdry X,H,J)$,
	is called the \emph{conformal infinity} of $g$.
\end{dfn}

Our fundamental questions on ACH metrics are the following. For a given $\overline{X}$,
does there exist an Einstein ACH metric on $X$ with prescribed conformal infinity?
If there does, how many are there essentially (i.e., up to the action of diffeomorphisms)?

Let us focus on the existence problem for the moment.
The Cheng--Yau theorem (Theorem \ref{thm:Cheng-Yau}) provides many examples of Einstein ACH metrics, but for
general infinity, only perturbative results are known.
Such results are given by Roth \cite{Roth-99-Thesis}, Biquard \cite{Biquard-00},
and the present author \cite{Matsumoto-preprint16}, which we shall now discuss.

In \cite{Roth-99-Thesis} and \cite{Biquard-00}, general perturbation theory is established.
Roth considered deformations of the Cheng--Yau metrics, while Biquard worked on those of
arbitrary Einstein ACH metrics.
It was shown that, in the both works, that if the given Einstein metric $g$ has negative sectional curvature everywhere,
then compatible almost CR structures nearby the conformal infinity of $g$ are also ``fillable'' with Einstein metrics.
More precisely, the following theorem holds.

\begin{thm}[Biquard \cite{Biquard-00}]
	\label{thm:Biquard-deformation}
	Let $g$ be an Einstein ACH metric on $X$, whose conformal infinity is denoted by $(\bdry X,H,J_0)$.
	Suppose that $g$ has negative sectional curvature.
	Then, if $\mathcal{J}$ is a sufficiently small $C^{2,\alpha}$ neighborhood of $J_0$ in the space of
	compatible almost CR structures on $H$,
	any $J\in\mathcal{J}$ is the conformal infinity of some Einstein ACH metric on $X$.
\end{thm}

In particular, Theorem \ref{thm:Biquard-deformation} is applicable to the complex hyperbolic metric
on the unit ball $B^n$ in $\mathbb{C}^n$ (note also that it is the Cheng--Yau metric of $B^n$).
Leaving the contact distribution $H$ unchanged is not an additional restriction,
because contact structures of closed manifolds are rigid.

Here is a very brief sketch of the construction (which is discussed more in the next section).
We first assign to each $J\in\mathcal{J}$ an approximate ACH solution $g_J$ of the Einstein equation
which satisfies
\begin{equation*}
	\Ric(g_J)+(n+1)g_J\in C^{0,\alpha}_\delta(X,S^2T^*X)
\end{equation*}
for some $\delta>0$ ($\delta$ must be independent of $J$).
We can do it in such a way that $g_J$ is smooth in $J$ and $g_{J_0}$ equals the original metric $g$.
Then we use functional analysis to show that, making $\mathcal{J}$ smaller if necessary,
for each $J\in\mathcal{J}$ one can find $\sigma\in C^{2,\alpha}_\delta(X,S^2T^*X)$ for which
$g'_J=g_J+\sigma$ satisfies $\Ric(g'_J)=-(n+1)g'_J$.
Since the modification term $\sigma$ belongs to $C^{2,\alpha}_\delta$,
$g'_J$ is still an ACH metric whose conformal infinity is $J$.

The negative curvature assumption is an easy sufficient condition that makes this plan work.
However, in practice, it is a nontrivial matter to check whether this condition is satisfied for a given $g$.
The following theorem shows that it is unnecessary for the Cheng--Yau metrics,
(at least) except for the two-dimensional case.

\begin{thm}[Matsumoto \cite{Matsumoto-preprint16}]
	\label{thm:main}
	Let $\Omega$ be a smoothly bounded strictly pseudoconvex domain in a Stein manifold of dimension $n\ge 3$,
	and $\mathcal{J}$ a sufficiently small $C^{2,\alpha}$ neighborhood of the induced CR structure $J_0$
	in the space of compatible almost CR structures on the natural contact distribution over $\bdry\Omega$.
	Then for each $J\in\mathcal{J}$, there is an Einstein ACH metric on $\Omega$ with conformal infinity $J$.
\end{thm}

There are such perturbation theorems also for Poincar\'e-Einstein metrics.
The possibility of deforming the real hyperbolic metric is shown by Graham--Lee \cite{Graham-Lee-91},
and in \cite{Biquard-00} it was pointed out that the negative curvature assumption is sufficient.
Lee \cite{Lee-06} showed that a weaker curvature assumption suffices
when the boundary conformal structure has nonnegative Yamabe constant.

In \cite{Biquard-00}, the local uniqueness of Einstein ACH metrics is also discussed.
By shrinking $\mathcal{J}$ if necessary,
the Einstein metric $g'_J$ constructed for each $J\in\mathcal{J}$ is the unique Einstein metric
modulo diffeomorphism action in a neighborhood of $g'_J$ in
$g'_J+C^{2,\alpha}_\delta(X,S^2T^*X)$ (for any $\delta>0$).
Probably the following refined question can be asked:
is there a neighborhood of $g$ in the unweighted H\"older space $C^{2,\alpha}(X,S^2T^*X)$ in which
there is only one Einstein metric for each conformal infinity?
To the author's knowledge, this is not yet settled so far.

\section{Ideas of the Proofs of Theorems \ref{thm:Biquard-deformation} and \ref{thm:main}}

The two theorems in the previous section are
reduced to the vanishing of the $L^2$ kernel of the ``linearized gauged Einstein operator''
acting on symmetric 2-tensors, which is
\begin{equation}
	\label{eq:linearized-gauged-Einstein-operator}
	P=\frac{1}{2}(\nabla^*\nabla-2\mathring{R}),
\end{equation}
where $g$ is the given Einstein ACH metric
and $\mathring{R}$ denotes the pointwise linear action of the curvature tensor of $g$.
Let us see how this reduction is carried out.

It is natural to study the linearization of the Einstein equation in order to deform Einstein metrics.
However, if we consider the Einstein equation itself, we encounter a difficulty
that originates from the diffeomorphism invariance of the equation.
One usually introduces an additional term to break this ``gauge invariance.''
Here we set, following \cite{Biquard-00},
\begin{equation*}
	\mathcal{E}_g(g'):=\Ric(g')+(n+1)g'+\delta_{g'}^*\mathcal{B}_g(g'),\qquad
	\mathcal{B}_g(g'):=\delta_gg'+\frac{1}{2}d\tr_gg'.
\end{equation*}
As long as we consider $g'$ in a small neighborhood of $g$ in $g+C^{2,\alpha}_\delta(X,S^2T^*X)$,
any solution of $\mathcal{E}_g(g')=0$ automatically satisfies $\mathcal{B}_g(g')=0$,
and hence it becomes an Einstein metric.

We apply the implicit function theorem to the mapping
\begin{equation*}
	\mathcal{J}\times C^{2,\alpha}_\delta(X,S^2T^*X)\to
	C^{0,\alpha}_\delta(X,S^2T^*X),\qquad
	(J,\sigma)\mapsto \mathcal{E}_{g_J}(g_J+\sigma)
\end{equation*}
at $(J_0,0)$, where $g_J$ is a family of approximate solutions as described in the previous section.
If the linearization of $\sigma\mapsto\mathcal{E}_g(g+\sigma)$ at $\sigma=0$,
which is the operator \eqref{eq:linearized-gauged-Einstein-operator}, is invertible,
then for each $J\in\mathcal{J}$ sufficiently close to $J_0$
there exists $\sigma\in C^{2,\alpha}_\delta(X,S^2T^*X)$ satisfying $\mathcal{E}_{g_J}(g_J+\sigma)=0$.
Thus it suffices to prove that \eqref{eq:linearized-gauged-Einstein-operator} is an isomorphism as the mapping
\begin{equation}
	\label{eq:linearized-gauged-Einstein-operator-Holder-spaces}
	P\colon C^{2,\alpha}_\delta(X,S^2T^*X)\to C^{0,\alpha}_\delta(X,S^2T^*X)
\end{equation}
for sufficiently small $\delta>0$.

An essential part is to show that
\eqref{eq:linearized-gauged-Einstein-operator-Holder-spaces} is an isomorphism for small $\delta>0$
if and only if
\begin{equation}
	\label{eq:linearized-gauged-Einstein-operator-L2-spaces}
	P\colon H^2(X,S^2T^*X)\to L^2(X,S^2T^*X)
\end{equation}
is isomorphic, where $H^2(X,S^2T^*X)$ denotes the $L^2$ Sobolev space of order $2$,
which is actually the domain of $P$ seen as an unbounded operator on $L^2(X,S^2T^*X)$.
It is easy to show that $P$ is a self-adjoint unbounded operator, and hence
\eqref{eq:linearized-gauged-Einstein-operator-L2-spaces} is isomorphic if the $L^2$ kernel vanishes.
The equivalence of \eqref{eq:linearized-gauged-Einstein-operator-Holder-spaces}
and \eqref{eq:linearized-gauged-Einstein-operator-L2-spaces} being isomorphic follows by
a certain parametrix construction, which makes good use of the geometry of ACH metrics, explained in \cite{Biquard-00}.
The exposition on the Poincar\'e-Einstein case in \cite{Lee-06} is also useful.

Consequently, it suffices to show that the $L^2$ kernel of $P$ is trivial.
When $g$ has negative sectional curvature, the vanishing can be proved by the following Bochner technique.
Note that any element of the kernel must be trace-free,
because if $\sigma=u g$ for some $u\in C^\infty(X)$ then $P\sigma=(\nabla^*\nabla u+2(n+1)u)g$,
and the operator $\nabla^*\nabla+2(n+1)$ acting on functions has trivial $L^2$ kernel.
Now if a general symmetric 2-tensor $\sigma$ is regarded as a 1-form with values in $T^*X$, then
$P$ satisfies the Weitzenb\"ock formula below given in terms of the exterior covariant differentiation $D$
(see \cite[12.69]{Besse-87}):
\begin{equation*}
	2P\sigma=(DD^*+D^*D)\sigma-\mathring{R}\sigma+(n+1)\sigma.
\end{equation*}
Moreover, there is also a pointwise estimate valid for trace-free $\sigma$ (\cite[12.71]{Besse-87}) that
\begin{equation*}
	\braket{\mathring{R}\sigma,\sigma}\le (n+1+2(n-1)K_\mathrm{max})\abs{\sigma}^2,
\end{equation*}
where $K_\mathrm{max}$ is the maximum of the sectional curvatures at a point.
Since the assumption implies that the sectional curvature is bounded from above by a negative constant
(by virtue of the asymptotic complex hyperbolicity),
one can deduce that the $L^2$ kernel of $P$ is trivial in this case.

When $g$ is the Cheng--Yau metric of a smoothly bounded strictly pseudoconvex domain $\Omega$ in a Stein manifold,
we argue as follows based on Koiso's observations \cite{Koiso-83} (see also Besse \cite[Section 12.J]{Besse-87}).
The K\"ahlerness of $g$ implies that $P$ respects the type decomposition of $\sigma$
into hermitian and anti-hermitian parts, and due to the Einstein condition,
$P$ on each type becomes a familiar operator.
If $\sigma$ is hermitian, then one may regard it as a $(1,1)$-form and we have
\begin{equation*}
	2P\sigma=(dd^*+d^*d)\sigma+2(n+1)\sigma.
\end{equation*}
This shows the vanishing of the hermitian part of the $L^2$ kernel.
On anti-hermitian symmetric 2-tensors, by regarding them as $(0,1)$-forms with values in
the holomorphic tangent bundle $T^{1,0}\Omega$, we obtain
\begin{equation*}
	2P\sigma=
	(\overline{\partial}\,\smash{\overline{\partial}}^*+\smash{\overline{\partial}}^*\overline{\partial})\sigma.
\end{equation*}
Therefore it suffices to show that there are no nontrivial
$L^2$ harmonic $(0,1)$-forms with values in $T^{1,0}\Omega$.
This allows one to restate the problem in terms of cohomology:
the isomorphicity of \eqref{eq:linearized-gauged-Einstein-operator-Holder-spaces} for small $\delta>0$ follows if
the $L^2$ Dolbeault cohomology $H^{0,1}_{(2)}(\Omega,T^{1,0}\Omega)$ vanishes.

It is a consequence of classical theory on Stein manifolds that the compactly supported cohomology
$H^{0,1}_c(\Omega,T^{1,0}\Omega)$ vanishes.
Moreover, it can be observed that the following sequence involving the inductive limit
$\varinjlim_K H^{0,1}_{(2)}(\Omega\setminus K;T^{1,0}(\Omega\setminus K))$,
where $K$ runs through all compact subsets of $\Omega$,
is exact:
\begin{equation*}
	\dotsb
	\to H_c^{0,1}(\Omega;T^{1,0}\Omega)
	\to H^{0,1}_{(2)}(\Omega;T^{1,0}\Omega)
	\to\varinjlim_K H^{0,1}_{(2)}(\Omega\setminus K;T^{1,0}(\Omega\setminus K))
	\to H_c^{0,2}(\Omega;T^{1,0}\Omega)
	\to \dotsb.
\end{equation*}
Therefore, $H^{0,1}_{(2)}(\Omega,T^{1,0}\Omega)=0$ follows if
\begin{equation}
	\label{eq:direct-limit}
	\varinjlim_K H^{0,1}_{(2)}(\Omega\setminus K,T^{1,0}(\Omega\setminus K))=0
\end{equation}
holds. We show \eqref{eq:direct-limit} by proving
\begin{equation}
	\label{eq:vanishing-in-collar-neighborhood}
	H^{0,1}_{(2)}(\mathcal{U},T^{1,0}\mathcal{U})=0,
\end{equation}
where $\mathcal{U}$ is
a sufficiently narrow collar neighborhood of $\bdry\Omega$ intersected with $\Omega$.
The vanishing \eqref{eq:vanishing-in-collar-neighborhood} is attacked by the usual technique of $L^2$ estimate,
but one needs to be careful because boundary integrals along the inner boundary of $\mathcal{U}$,
which is strictly pseudoconcave, comes into play.
The $L^2$ estimate so obtained is actually
sufficient to prove \eqref{eq:vanishing-in-collar-neighborhood} only when $n\ge 4$.
When $n=3$, one needs to work with a weighted $L^2$ cohomology instead.

\section{Problems}

An obvious problem related to Theorem \ref{thm:main} is to clarify what happens in the two-dimensional case.
The author expects (perhaps optimistically) that finally
one can simply remove the assumption $n\ge 3$ from the theorem,
for it is at least true for the unit ball by Theorem \ref{thm:Biquard-deformation},
and there seems to be no reason to expect that it fails for general strictly pseudoconvex domains.

A more challenging issue about existence is how we can construct Einstein ACH metrics
for compatible almost CR structures which are far from those that are known to be ``fillable.''
The corresponding problem in the Poincar\'e-Einstein setting is also a long-standing one.
One should be aware that there is a recent \emph{nonexistence} result
by Gursky--Han \cite{Gursky-Han-17} in the latter setting.

Turning to the uniqueness of Einstein fillings for a given conformal infinity,
in the Poincar\'e-Einstein case,
an example of Hawking--Page \cite{Hawking-Page-82} exhibits that it fails in general
(see Anderson \cite{Anderson-03} for further explanation).
A similar nonuniqueness example for ACH-Einstein metrics will be of great interest,
as well as uniqueness results under some assumption.
There is also a room for further investigations
about local uniqueness as mentioned at the end of Section \ref{sec:ACH-metrics}.

A typical application of Poincar\'e-Einstein metrics is the construction of conformally invariant objects on the
boundary, and there is a similar story for ACH-Einstein metrics.
For this purpose the determination of the asymptotic behavior of the metric in terms of
the boundary geometry is important, and its formal aspects are studied by
Fefferman--Graham \cites{Fefferman-Graham-85,Fefferman-Graham-12} for the Poincar\'e-Einstein metrics,
by Fefferman \cite{Fefferman-76} and Graham \cite{Graham-87} for the Cheng--Yau metrics
(as mentioned in Section \ref{sec:Introduction}),
and by Biquard--Herzlich \cite{Biquard-Herzlich-05} and the author \cite{Matsumoto-14}
(see also \cite{Matsumoto-13-Thesis}) for general ACH-Einstein metrics.
There is a tremendous amount of literature regarding constructions of conformal invariants based on
\cites{Fefferman-Graham-85,Fefferman-Graham-12},
while in the CR case such constructions are discussed in, e.g.,
\cites{Fefferman-79,Graham-87,Graham-87-Invariants,Burns-Epstein-90,Bailey-Eastwood-Graham-94,Hirachi-00,Fefferman-Hirachi-03,Biquard-Herzlich-05,Gover-Graham-05,Hislop-Perry-Tang-08,Case-Yang-13,Hirachi-14,Marugame-16,Matsumoto-16,Hirachi-Marugame-Matsumoto-17,Takeuchi-17,Case-Gover-17preprint,Marugame-18,Takeuchi-18,Takeuchi-18toappear,Marugame-18preprint}. Further developments along this line are anticipated.
It would be also very interesting if there is some invariant construction that needs
global considerations on Einstein metrics in an essential way.

Now let us take notice of the fact
that the Cheng--Yau metrics come with complex structures with respect to which they are K\"ahler.
As a problem without any counterpart in the Poincar\'e-Einstein setting,
it may be interesting to look for a canonical way to determine a good almost complex structure on
a manifold equipped with an ACH-Einstein metric.
That is to say, Riemannian metrics may not be the ``best'' filling geometric structure inside CR manifolds.
It seems to the author that this idea is backed up by
the fact that the Einstein deformation problem is recast in the proof of Theorem \ref{thm:main}
in terms of harmonic $(0,1)$-forms with values in the holomorphic tangent bundle.

Finally, the author would like to remark once again that it should also be fruitful to examine
geometries modelled on other symmetric spaces.
Interested readers are referred to Biquard \cites{Biquard-99,Biquard-00},
Biquard--Mazzeo \cite{Biquard-Mazzeo-06,Biquard-Mazzeo-11}, and references therein.

\section*{Acknowledgments}

I wish to express my gratitude to the hospitality of Stanford University, where I was working as a visiting member
when the manuscript was written, and I am deeply grateful to Rafe Mazzeo for hosting the visit, for discussions,
and for encouragements.
I would also appreciate the careful reading of the manuscript by the reviewer.
This work was partially supported by JSPS KAKENHI Grant Number JP17K14189 and
JSPS Overseas Research Fellowship.

\bibliography{myrefs}

\end{document}